\documentclass[a4paper,11pt]{article}
\usepackage{amssymb,latexsym,epsf,graphicx}

\newtheorem{thm}{Theorem}[section]
\newtheorem{lem}[thm]{Lemma}

\newtheorem{cor}[thm]{Corollary}
\newtheorem{prop}[thm]{Proposition}

\newtheorem{rem}[thm]{Remark}

\title{Chebyshev polynomials, quadratic surds
and a variation of Pascal's triangle}
\author{Roland Bacher\footnote{This work has been partially supported by the LabEx PERSYVAL-Lab (ANR--11-LABX-0025). The author is a member of the project-team GALOIS supported by this LabEx.}
}

\begin{document}
\maketitle

{\sl Abstract\footnote{Keywords: Chebyshev polynomial, continued fraction, binomial coefficient.
Math. class: Primary: 33C45, Secondary: 11A55, 11B65}: Using iterated Chebyshev polynomials
of both kinds, we construct rational fractions which are convergents 
of the smallest root of $x^2-\alpha x+1$ for $\alpha=3,4,5,\dots$.

Some of the underlying identities suggest an identity involving binomial
coefficients which leads to a triangular array sharing many properties
with Pascal's triangle.}



\section{Introduction}

refs: Watson-Whittacker, qq chose sur fractions continues, polys de Cheb.

Chebyshev polynomials of the first kind $T_0,T_1,\dots$
and of the second kind $U_0,U_1,\dots$ have recursive definitions given by
$$T_0(x)=1,T_1(x)=x,T_{n+1}(x)=2xT_n(x)-T_{n-1}(x),n\geq 1,$$
and by
$$U_0(x)=1,U_1(x)=2x,U_{n+1}(x)=2xU_n(x)-U_{n-1}(x),n\geq 1.$$

We write in the sequel always simply $T_n,U_n$ for $T_n(x),U_n(x)$.
The polynomial sequences $T_0,T_1,\dots$ and $U_0,U_1,\dots$ satisfy 
the same linear recursion relation with characteristic polynomial 
$Z^2-2xZ+1$. This implies easily the formulae
\begin{eqnarray}\label{formusqrrtT}
T_n&=&\frac{\left(x+\sqrt{x^2-1}\right)^n+\left(x-\sqrt{x^2-1}\right)^n}{2}
\end{eqnarray}
and 
\begin{eqnarray}\label{formusqrrtU}
U_n&=&\frac{\left(x+\sqrt{x^2-1}\right)^{n+1}-\left(x-\sqrt{x^2-1}\right)^{n+1}}{2\sqrt{x^2-1}}.
\end{eqnarray}

The identity $2\cos x\cos nx=\cos(n+1)x+\cos(n-1)x$ (together with the 
initial values $T_0=1,T_1=x$) implies by an easy induction 
$T_n(\cos x)=\cos nx$ (which is often used for defining 
Chebyshev polynomials of the first kind). As a consequence,
we have the identity
\begin{eqnarray}\label{idTnmeqTnoTm}
T_{nm}&=&T_n\circ T_m.
\end{eqnarray}


For $n\in \mathbb N$ and $d\geq 1$ we introduce the rational fraction
\begin{eqnarray}\label{formulaSn}
S_{n,d}=S_{n,d}(x)&=&\sum_{k=0}^nU_{d-1}\circ T_{(d+1)^k}\prod_{j=0}^k
\frac{1}{U_d\circ T_{(d+1)^j}}.
\end{eqnarray}

\begin{thm}\label{thmmain} We have for all $n\in\mathbb N$ and for all $d\geq 1$ the 
identity
\begin{eqnarray}\label{thmmaineq}
S_{n,d}^2-2xS_{n,d}+1&=&\left(\prod_{j=0}^n
\frac{1}{U_d\circ T_{(d+1)^j}}\right)^2.
\end{eqnarray}
\end{thm}

\begin{cor}\label{corlimval} For $x$ evaluated to a real number in $\mathbb R\setminus
[-1,1]$, the limit $\lim_{n+d\rightarrow \infty}S_{n,d}$ exists and is 
given by the root $\hbox{sign}(x)\left(\vert x\vert-\sqrt{x^2-1}\right)$
(where $\hbox{sign}(x)\in\{\pm 1\}$ denotes the sign of $x$)
closest to zero of $X^2-2x X+1$.
\end{cor}

The following result expresses $S_{n,d}$ as a simple fraction:
\begin{thm}\label{thmb} We have
$$S_{n,d}=\frac{U_{(d+1)^{n+1}-2}}{U_{(d+1)^{n+1}-1}}.$$
\end{thm}

Corollary \ref{corlimval} is now an almost immediate consequence of 
Formula (\ref{formusqrrtU}) and Theorem \ref{thmb}.

Note that Formula (\ref{formulaSn}), perhaps computed 
using iteratively (\ref{idTnmeqTnoTm}), is perhaps better suited for 
computations than the simpler expression given 
by Theorem \ref{thmb}.

Finally we have the following result:
\begin{thm}\label{thmc}
We have 
$$\frac{U_n}{U_{n+1}}=\left[0;2x-1,1\left(,2(x-1),1\right)^n\right].$$
\end{thm}

Theorem \ref{thmb} and Theorem \ref{thmc} together imply
that the evaluation at an integer $x\geq 2$ of
$S_{n,d}$ is a convergent of the real quadratic surd
$$x+\sqrt{x^2-1}=[0;2x-1,1,2(x-1),1,2(x-1)\overline{,1,2(x-1)}]$$
with minimal polynomial $X-2xX+1\in\mathbb Z[X]$.

In a last only losely related part we study some identities involving 
binomial coefficients (obtained by expressing Chebyshev polynomials in terms
of binomials coefficients and generalizing some of the previously obtained
identities). They lead to the array of numbers 
$$\begin{array}{cccccccccccccccccccc}
&&&&&&&&&1\\
&&&&&&&&-1&&1\\
&&&&&&&1&&0&&1\\
&&&&&&-1&&1&&1&&1\\
&&&&&1&&0&&2&&2&&1\\
&&&&-1&&1&&2&&4&&3&&1\\
&&&1&&0&&3&&6&&7&&4&&1\\
&&-1&&1&&3&&9&&13&&11&&5&&1\\
&1&&0&&4&&12&&22&&24&&16&&6&&1\\
\end{array}$$
mimicking several aspects of Pascal's triangle.

The sequel of the paper is organised as follows:

Section \ref{sectVajda} describes and proves useful identities among Chebyshev polynomials.

Section \ref{sectmainb} contains a proof of Theorem \ref{thmmain}, \ref{thmb}
and \ref{thmc}.

A final Section describes a few relations with binomial coefficients
and studies a few features of the above analogue of Pascal's triangle.

\section{Analogues of Vajda's identity for Chebyshev polynomials}
\label{sectVajda}

Fibonacci numbers $F_1=F_2=1,F_{n+1}=F_n+F_{n-1}$ satisfy Vajda's identity 
$$F_{n+i}F_{n+j}-F_nF_{n+i+j}=(-1)^nF_iF_j.$$

The following result describes analogues for Chebyshev polynomials:


\begin{thm}\label{thmChebid} We have the following identities for Chebyshev polynomials:
\begin{eqnarray*}
U_{n+i}U_{n+j}-U_{n-1}U_{n+1+i+j}&=&U_iU_j\\
T_{n+i}T_{n+j}-T_{n-1}T_{n+1+i+j}&=&(1-x^2)U_iU_j\\
T_{n+i}U_{n+j}-U_{n-1}T_{n+1+i+j}&=&T_iU_j\\
T_{n+i}U_{n+j}-T_{n-1}U_{n+1+i+j}&=&-U_iT_{j+2}\\
T_{n+i}T_{n+j}-(x^2-1)U_{n-1}U_{n-1+i+j}&=&T_iT_j\\
\end{eqnarray*}
\end{thm}

The case $i=j=0$ of the first identity
specialises to the so-called Cassini-Simpson
identity $U_n^2-U_{n+1}U_{n-1}=1$ and implies Turan's 
inequality $U_n^2(x)>U_{n+1}(x)U_{n-1}(x)$ for all real $x$.

The last equation generalises the instance 
$T_n^2-(x^2-1)U_{n-1}^2=1$ (corresponding to $i=j=0$) of Pell's equation.

Only a few cases of the first identity will in fact be used in the sequel.

{\bf Proof of Theorem \ref{thmChebid}}
We consider
$$R(n,i,j)=U_{n+i}U_{n+j}-U_{n-1}U_{n+1+i+j}-U_iU_j.$$
We have to show that $R(n,i,j)=0$ for all $n\geq 1$ and for all 
$i,j\in \mathbb N$. Using the recursion relation $U_{n+1}=2xU_n-U_{n-1}$
on all terms depending on $i$, respectively depending on $j$,
we see that it is enough to prove the equalities $R(n,i,j)=0$
for $i,j\in\{0,1\}$. Using the obvious identity $R(n,i,j)=R(n,j,i)$
we are left with three cases: $R(n,0,0),R(n,0,1)$ and $R(n,1,1)$.

The computation
\begin{eqnarray*}
&&U_n^2-U_{n-1}U_{n+1}-U_0^2\\
&=&U_n(2xU_{n-1}-U_{n-2})-U_{n-1}(2xU_n-U_{n-1})-1\\
&=&U_{n-1}^2-U_nU_{n-2}-1
\end{eqnarray*}
shows $R(n,0,0)=R(n-1,0,0)$ for $n\geq 2$.
Similarly,
\begin{eqnarray*}
&&U_{n+1}U_n-U_{n-1}U_{n+2}-U_1U_0\\
&=&(2xU_n-U_{n-1})U_n-U_{n-1}(2xU_{n+1}-U_n)-2x\\
&=&2x(U_n^2-U_{n-1}U_{n+1}-1)\\
\end{eqnarray*}
shows $R(n,1,0)=R(n,0,1)=2xR(n,0,0)$.
Finally, the identities
\begin{eqnarray*}
&&U_{n+1}U_{n+1}-U_{n-1}U_{n+3}-U_1^2\\
&=&U_{n+1}(2xU_n-U_{n-1})-U_{n-1}(2xU_{n+2}-U_{n+1})-4x^2\\
&=&2x(U_{n+1}U_n-U_{n-1}U_{n+2}-2x)
\end{eqnarray*}
show $R(n,1,1)=2xR(n,1,0)$.
It is now enough to check that $R(n,i,j)=0$ for $n\in\{1,2\}$ and 
$i,j\in\{0,1\}$. 

Proofs of the remaining identities are 
similar.\hfill$\Box$

\begin{rem} Short direct proofs of Theorem \ref{thmChebid} can be obtained using Formulae (\ref{formusqrrtT})
and (\ref{formusqrrtU}).
\end{rem}

\section{Proof of Theorem \ref{thmmain}, \ref{thmb}
and  \ref{thmc}}\label{sectmainb}

\subsection{Useful identities}
The following result is well-known:
\begin{lem}\label{lemUTid} We have for all $d\geq 1$ the identity
$$xU_d=U_{d-1}+T_{d+1}.$$
\end{lem}

The proof is an easy induction left to the reader.\hfill$\Box$

\begin{lem}\label{lemUUid} We have for all $d\geq 1$ the identity
$$U_d^2=U_{d-1}^2+2T_{d+1}U_{d-1}+1.$$
\end{lem}

\noindent{\bf Proof} Using Lemma \ref{lemUTid} and the recursive 
definition of $U$, Lemma \ref{lemUUid} is equivalent to 
\begin{eqnarray*}
U_d^2&=&U_{d-1}^2+2(xU_d-U_{d-1})U_{d-1}+1\\
&=&-U_{d-1}^2+(2xU_d)U_{d-1}+1\\
&=&-U_{d-1}^2+(U_{d+1}+U_{d-1})U_{d-1}+1\\
&=&U_{d+1}U_{d-1}+1
\end{eqnarray*}
which is a special case of the 
first equality in Theorem \ref{thmChebid}.\hfill$\Box$

\begin{lem}\label{lem2TU} We have
$$2T_dU_n=U_{n+d}+U_{n-d}$$
for all $n\in\mathbb N$ and for all $d\in\{0,1,\dots,n\}$.
\end{lem}

\noindent{\bf Proof} We set 
$$R(n,d)=2T_dU_n-U_{n+d}-U_{n-d}.$$
Since $T_0=1$ and $T_1=x$ we have $R(n,0)=R(n,1)=0$.
The identities
\begin{eqnarray*}
R(n,d)&=&2T_dU_n-U_{n+d}-U_{n-d}\\
&=&2(2xT_{d-1}-T_{d-2})U_n-(2xU_{n+d-1}-U_{n+d-2})-(2xU_{n-d+1}-U_{n-d+2})\\
&=&2xR(n,d-1)-R(n,d-2)
\end{eqnarray*}
finish the proof.\hfill$\Box$

\begin{lem}\label{lemUUT} We have for all $n\geq 1$ the identity
$$U_{2n-1}=2T_nU_{n-1}.$$
\end{lem}

\noindent{\bf Proof} Equality holds for $n=1$. Using Lemma (\ref{lemUTid}) we have for $n\geq 2$ the identities
\begin{eqnarray*}
&&U_{2n-1}-2T_nU_{n-1}\\
&=&U_{2n-1}-2(xU_{n-1}-U_{n-2})U_{n-1}\\
&=&U_{2n-1}-(2xU_{n-1}-U_{n-2})U_{n-1}+U_{n-1}U_{n-2}\\
&=&-U_nU_{n-1}+U_{2n-1}+U_{n-1}U_{n-2}\\
&=&-(U_{1+(n-1)}U_{1+(n-2)}-U_{1-1}U_{1+1+(n-1)+(n-2)}-U_{n-1}U_{n-2}).
\end{eqnarray*}
The last expression equals zero by the first identity of Theorem \ref{thmChebid}.
\hfill$\Box$

\begin{lem}\label{lemUUTUUT} We have
$$U_{(n-1)d-1}U_{n-1}\circ T_d=U_{nd-1}U_{n-2}\circ T_d$$
for all $n\geq 2$ and for all $d\geq 1$.
\end{lem}

\noindent{\bf Proof} The case $n=2$ boils down $U_{2d-1}=2U_{d-1}T_d$ 
which holds by Lemma (\ref{lemUUT}).

Adding to
$$0=(2T_dU_{(n-1)d-1}-U_{(n-2)d-1}-U_{nd-1})U_{n-2}\circ T_d$$
which holds by Lemma (\ref{lem2TU}) the induction hypothesis 
we get
\begin{eqnarray*}
0&=&(2T_dU_{(n-1)d-1}-U_{(n-2)d-1}-U_{nd-1})U_{n-2}\circ T_d\\
&&\quad +U_{(n-2)d-1}U_{n-2}\circ T_d-U_{(n-1)d-1}U_{n-3}\circ T_d\\
&=&U_{(n-1)d-1}\left(2T_dU_{n-2}\circ T_d-U_{n-3}\circ T_d\right)-U_{nd-1}U_{n-2}\circ T_d\\
&=&U_{(n-1)d-1}U_{n-1}\circ T_d-U_{nd-1}U_{n-2}\circ T_d\\
\end{eqnarray*}
which ends the proof.\hfill$\Box$


\subsection{Proof of Theorem \ref{thmmain}}

We prove first that equation (\ref{thmmaineq}) holds for $n=0$.
Multiplying the left-side of equation (\ref{thmmaineq})
by $U_d^2$, we get
\begin{eqnarray*}
&&U_{d-1}^2-2xU_{d-1}U_d+U_d^2\\
&=&U_{d-1}^2-U_{d-1}(U_{d+1}+U_{d-1})+U_d^2\\
&=&U_d^2-U_{d-1}U_{d+1}\\
&=&U_0^2=1
\end{eqnarray*}
by applying the recursive definition of $U_i$ and the first 
identity of Theorem \ref{thmChebid} with $n=d,i=j=0$.

Setting $x=T_{d+1}$ in equation (\ref{thmmaineq}) and dividing the 
result by $U_d^2$, we have now by induction
\begin{eqnarray}\label{eqnindstep}
\left(\frac{S_{n,d}\circ T_{d+1}}{U_d}\right)^2-2T_{d+1}\frac{S_{n,d}\circ T_{d+1}}{U_d^2}+\frac{1}{U_{d}^2}&=&\left(\prod_{j=0}^{n+1}\frac{1}{U_d\circ T_{(d+1)^j}}\right)^2
\end{eqnarray}
where we have also used (\ref{idTnmeqTnoTm}) on the right side.
We rewrite now the obvious identity
\begin{eqnarray}\label{idinductSn}
S_{n+1,d}&=&\frac{U_{d-1}+S_{n,d}\circ T_{d+1}}{U_d}.
\end{eqnarray}
as
\begin{eqnarray}\label{idinducSnb}
S_{n,d}\circ T_{d+1}=U_dS_{n+1,d}-U_{d-1}.
\end{eqnarray}
Using (\ref{idinducSnb}) the left side of (\ref{eqnindstep}) equals
\begin{eqnarray*}
&&\left(\frac{U_dS_{n+1,d}-U_{d-1}}{U_d}\right)^2-2T_{d+1}\frac{U_dS_{n+1,d}-U_{d-1}}{U_d^2}+\frac{1}{U_d^2}\\
&=&S_{n+1,d}^2-2\frac{U_{d-1}}{U_d}S_{n+1,d}+\frac{U_{d-1}^2}{U_d^2}-
2\frac{T_{d+1}}{U_d}S_{n+1,d}+2T_{d+1}\frac{U_{d-1}}{U_d^2}+\frac{1}{U_d^2}
\end{eqnarray*}
Since 
$$-2\frac{U_{d-1}}{U_d}S_{n+1,d}-2\frac{T_{d+1}}{U_d}S_{n+1,d}
=-2xS_{n+1,d}$$
by Lemma (\ref{lemUTid}) and
$$\frac{U_{d-1}^2}{U_d^2}+2T_{d+1}\frac{U_{d-1}}{U_d^2}+\frac{1}{U_d^2}
=1$$
by Lemma (\ref{lemUUid}), we get
finally $S_{n+1,d}^2-2xS_{n+1,d}+1$ for the left side of (\ref{eqnindstep}).
This ends the proof.\hfill$\Box$

\subsection{Proof of Theorem \ref{thmb}}
Equality holds obviously for $n=0$. Applying the induction hypothesis to
(\ref{idinductSn})
we have to establish the equality
$$U_d\frac{U_{(d+1)^{n+2}-2}}{U_{(d+1)^{n+2}-1}}=U_{d-1}+
\frac{U_{(d+1)^{n+1}-2}\circ T_{d+1}}{U_{(d+1)^{n+1}-1}\circ T_{d+1}}$$
equivalent to 
\begin{eqnarray*}
0&=&(U_dU_{(d+1)^{n+2}-2}-U_{d-1}U_{(d+1)^{n+2}-1})U_{(d+1)^{n+1}-1}\circ T_{d+1}\\
&&\ -U_{(d+1)^{n+2}-1}U_{(d+1)^{n+1}-2}\circ T_{d+1}\\
&=&U_{(d+1)^{n+2}-2-d}U_{(d+1)^{n+1}-1}\circ T_{d+1}\\
&&\ -U_{(d+1)^{n+2}-1}U_{(d+1)^{n+1}-2}\circ T_{d+1}\\
\end{eqnarray*}
where we have applied the first identity of Theorem \ref{thmChebid} 
with $n=d,i=0,j=(d+1)^{n+2}-2-d$.
The identity
\begin{eqnarray*}
0&=&U_{(d+1)^{n+2}-2-d}U_{(d+1)^{n+1}-1}\circ T_{d+1}\\
&&\ -U_{(d+1)^{n+2}-1}U_{(d+1)^{n+1}-2}\circ T_{d+1}\\
\end{eqnarray*}
is now the case $(n,d)=((d-1)^{n+1},d+1)$ of 
Lemma \ref{lemUUTUUT}.\hfill$\Box$

\subsection{Proof of Theorem \ref{thmc}}

\noindent{\bf Proof of Theorem \ref{thmc}} Equality holds for $n=0$.

Setting $\gamma_n=\left[0;2x-1,1\left(,2(x-1),1\right)^n\right]$ we have
$$\frac{\gamma_n}{1-\gamma_n}=\left[0;2x-2,1\left(,2(x-1),1\right)^{n-1}
\right]$$
showing 
\begin{eqnarray*}
\gamma_{n+1}&=&\frac{1}{2x-1+\frac{1}{1+\frac{\gamma_n}{1-\gamma_n}}}\\
&=&\frac{1}{2x-\gamma_n}.
\end{eqnarray*}
The result follows now by induction from the trivial identities
$$\frac{1}{2x-\frac{U_n}{U_{n+1}}}=\frac{U_{n+1}}{2x U_{n+1}-U_n}=
\frac{U_{n+1}}{U_{n+2}}.$$
\hfill$\Box$

An easy computation shows 
the continued fraction expansion
$$x-\sqrt{x^2-1}=[0;2x-1,1,2(x-1),1,2(x-1)\overline{,1,2(x-1)}].$$
for $x\in\{2,3,4,\dots\}$. Equality follows thus from analytic
continuation whenever both sides make sense.

Combining Theorem \ref{thmb} and Theorem \ref{thmc} we see that 
$$S_{n,d}=\left[0;2x-1,1\left(,2(x-1),1\right)^{(d+1)^{n+1}-2}\right]$$
(using a hopefully self-explanatory notation)
is a convergent of $x-\sqrt{x^2-1}$ for $x=2,3,\dots$.

\section{A sum of products of two binomial coefficients}

\subsection{Coefficients of Chebyshev polynomials}

\begin{lem}\label{lemUnid} Explicit expressions for 
coefficients of Chebyshev polynomials are given by the formulae
\begin{eqnarray*}
T_n&=&\frac{1}{2}\sum_{k=0}^{\lfloor n/2\rfloor}
(-1)^k\left({n+1-k\choose k}-{n-1-k\choose k-2}\right)(2x)^{n-2k},\\
U_n&=&\sum_{k=0}^{\lfloor n/2\rfloor}
(-1)^k{n-k\choose k}(2x)^{n-2k}
\end{eqnarray*}
(using the conventions ${-1\choose -2}=-1$, ${\not=-1,-2\choose -2}=0$,
${-1\choose -1}=1$, ${\not=-1\choose -1}=0$).
\end{lem}

\noindent{\bf Proof} The formulae hold obviously for $T_0,T_1$ and 
$U_0,U_1$. We have now
\begin{eqnarray*}
T_{n+1}&=&2xT_n-T_{n-1}\\
&=&2x\frac{1}{2}\sum_{k=0}^{\lfloor n/2\rfloor}
(-1)^k\left({n+1-k\choose k}-{n-1-k\choose k-2}\right)(2x)^{n-2k}\\
&&\quad -\frac{1}{2}\sum_{k=0}^{\lfloor (n-1)/2\rfloor}
(-1)^k\left({n-k\choose k}-{n-2-k\choose k-2}\right)(2x)^{n-1-2k}\\
&=&\frac{1}{2}\sum_k
(-1)^k\left({n+1-k\choose k}+{n+1-k\choose k-1}\right)(2x)^{n+1-2k}\\
&&\quad -\frac{1}{2}\sum_k
(-1)^k\left({n-1-k\choose k-2}+{n-1-k\choose k-3}\right)(2x)^{n+1-2k}\\
&=&\frac{1}{2}\sum_k
(-1)^k\left({n+2-k\choose k}-{n-k\choose k-2}\right)(2x)^{n+1-2k}
\end{eqnarray*}
and
\begin{eqnarray*}
U_{n+1}&=&2xU_n-U_{n-1}\\
&=&2x\sum_{k=0}^{\lfloor n/2\rfloor}
(-1)^k{n-k\choose k}(2x)^{n-2k}\\
&&\quad -\sum_{k=0}^{\lfloor (n-1)/2\rfloor}
(-1)^k{n-1-k\choose k}(2x)^{n-1-2k}\\
&=&\sum_k(-1)^k\left({n-k\choose k}+{n-k\choose k-1}\right)(2x)^{n+1-2k}\\
&=&\sum_k(-1)^k{n+1-k\choose k}(2x)^{n+1-2k}.
\end{eqnarray*}
These identities imply the result by induction.\hfill$\Box$

\subsection{A curious identity}

Rewriting Chebyshev polynomials in terms of binomial coefficients 
using the identities of Lemma \ref{lemUnid},
some identities among Chebyshev polynomials are special 
cases of the following result.

\begin{thm}\label{thmidbin} The expression 
$$f(a,d,n)_x=\sum_{k=0}^{d-n}{a+d+x-k\choose k}{d+k-x\choose d-n-k}$$
is constant in $x$ and depends only on $a,d\in\mathbb C$ and $n\in
d-\mathbb N=\{d,d-1,d-2,d-3,\dots\}$.
\end{thm}

Observe that all values $f(a,d,n)_x$ are determined by
the values $f(0,d,n)_x$ using the trivial identity
\begin{eqnarray}\label{idftriv}
f(a,d,n)_x&=&f(a-2c,d+c,n+c)_{x+c}
\end{eqnarray}
with $c=a/2$.

Theorem \ref{thmidbin} implies that 
$Q_N(z)=f(0,z/2,z/2-N)_*$ is a polynomial in $\mathbb Q[z]$
of degree $N$ such that $Q_N(\mathbb Z)\subset \mathbb Z$.

\begin{lem} We have the identities
\begin{eqnarray}\label{idfadn1}
f(a,d,n)_x&=&f(a-1,d,n)_x+f(a-1,d,n+1)_{x-1}
\end{eqnarray}
and 
\begin{eqnarray}\label{idfadn2}
f(a,d,n)_x&=&f(a-1,d,n)_{x+1}+f(a-1,d,n+1)_{x+1}
\end{eqnarray}
\end{lem}

\noindent{\bf Proof} Follows from the computations
\begin{eqnarray*}
&&f(a,d,n)_x\\
&=&\sum_{k=0}^{d-n}\left({a-1+d+x-k\choose k}+{a-1+d+x-k\choose k-1}\right)
{d+k-x\choose d-n-k}\\
&=&f(a-1,d,n)_x\\
&&\quad +\sum_k{a-1+d+x-1-(k-1)\choose k-1}
{d+(k-1)-(x-1)\choose d-(n+1)-(k-1)}\\
&=&f(a-1,d,n)_x+f(a-1,d,n+1)_{x-1}
\end{eqnarray*}
and
\begin{eqnarray*}
&&f(a,d,n)_x\\
&=&\sum_{k=0}^{d-n}{a+d+x-k\choose k}\left(
{d+k-x-1\choose d-n-k}+{d+k-x-1\choose d-n-1-k}\right)\\
&=&f(a-1,d,n)_{x+1}+f(a-1,d,n+1)_{x+1}
\end{eqnarray*}
\hfill$\Box$

\noindent{\bf Proof of Theorem \ref{thmidbin}} Since
${x\choose k}=\frac{x(x-1)\cdots (x-k+1)}{k!}$,
the function $f(a,d,n)_x$ is a polynomial of degree at most 
$d-n$ in $x$. It is thus independent of $x$ for $n=d$.
Subtracting equation (\ref{idfadn1}) from (\ref{idfadn2})
we get 
$$f(a-1,d,n)_{x+1}-f(a-1,d,n)_x=f(a-1,d,n+1)_{x-1}-f(a-1,d,n+1)_{x+1}$$
which implies the result by induction on $d-n$.\hfill$\Box$

\subsection{A few properties of $f(0,d,n)$}

The numbers 
\begin{eqnarray*}
l_{i,j}&=&f(0,(i-1)/2,j-(i+1)/2)\\
&=&\sum_{k=0}^{i-j}{\frac{i-1}{2}+x-k\choose k}{\frac{i-1}{2}+k-x\choose i-j-k}
\end{eqnarray*} 
with $i\in \mathbb N$ and 
$j\in\{0,\dots,i\}$ (and $x$ arbitrary)
form the ``Pascal-like'' triangle:
$$\begin{array}{cccccccccccccccccccc}
&&&&&&&&&1\\
&&&&&&&&-1&&1\\
&&&&&&&1&&0&&1\\
&&&&&&-1&&1&&1&&1\\
&&&&&1&&0&&2&&2&&1\\
&&&&-1&&1&&2&&4&&3&&1\\
&&&1&&0&&3&&6&&7&&4&&1\\
&&-1&&1&&3&&9&&13&&11&&5&&1\\
&1&&0&&4&&12&&22&&24&&16&&6&&1\\
\end{array}$$
as shown by the following result:

\begin{prop}\label{propladd} We have $l_{i,0}=(-1)^i,\ l_{i,i}=1$ and 
$l_{i,j}=l_{i-1,j-1}+l_{i-1,j}$.
\end{prop}

\noindent{\bf Proof}
Using $x=\frac{i-1}{2}$ we get for $j=0$ the evaluation
\begin{eqnarray*}
f\left(0,\frac{i-1}{2},-\frac{i+1}{2}\right)&=&
\sum_{k=0}^i{i-1-k\choose k}{k\choose i-k}\\
&=&{-1\choose i}{i\choose 0}=(-1)^i.
\end{eqnarray*}
For $i=j$ we get $l_{i,i}=f\left(0,\frac{i-1}{2},\frac{i-1}{2}-\frac{i+1}{2}\right)=1$
by a trivial computation.

Using (\ref{idfadn1}) followed by two applications 
of (\ref{idftriv}) with $c=-\frac{1}{2}$ we have
\begin{eqnarray*}
l_{i,j}&=&f\left(0,\frac{i-1}{2},j-\frac{i+1}{2}\right)\\
&=&f\left(-1,\frac{i-1}{2},j-\frac{i+1}{2}\right)+
f\left(-1,\frac{i-1}{2},j-\frac{i+1}{2}+1\right)\\
&=&f\left(0,\frac{i-1}{2}-\frac{1}{2},j-\frac{i+1}{2}-\frac{1}{2}\right)+
f\left(0,\frac{i-1}{2}-\frac{1}{2},j-\frac{i+1}{2}+\frac{1}{2}\right)\\
&=&f\left(0,\frac{(i-1)-1}{2},j-1-\frac{(i-1)+1}{2}\right)\\
&&\quad +
f\left(0,\frac{(i-1)-1}{2},j-\frac{(i-1)+1}{2}\right)\\
&=&l_{i-1,j-1}+l_{i-1,j}
\end{eqnarray*}
which proves the result.\hfill$\Box$

\subsection{An $LU$-decomposition}

Interpreting the integers $l_{i,j}, i,j\geq 0$ as the coefficients of an infinite
unipotent matrix $L$ and introducing similarly the matrix $M$
with coefficients $M_{i,j}=l_{i+j,j}$, we have the following result:

\begin{prop}\label{propLU} We have 
$$M=LU$$
where $U$ is the upper-triangular matrix with coefficients $U_{i,j}={j\choose i},
i,j\geq 0$.

In particular, we have $\det(M(n))=1$ where $M(n)$ is the square matrix
consisting of the first $n$ rows and columns of $M$.
\end{prop}

Proposition \ref{propLU} is a special case of the following more general 
result:

We associate two infinite matrices to an infinite sequence 
$\alpha_0=1,\alpha_1,\alpha_2,\dots$ in a commutative ring with $1$ as follows:

The first matrix $M(\alpha)$ with coefficients $M_{i,j}$ indexed by 
$i,j\in \mathbb N$ is defined recursively by 
$$M_{0,0}=1,M_{0,j}=1,M_{i,0}=\alpha_i,M_{i,j}=M_{i-1,j}+M_{i,j-1}, i,j\geq 1.$$
The coefficients $M_{i,j}$ for $j>0$ are also given by the formula
$$M_{i,j}=\sum_{k=0}^i{k+j-1\choose k}\alpha_{i-k}.$$
The second matrix is the unipotent lower-triangular matrix $L(\alpha)$
with lower triangular coefficients $L_{i,j}=M_{i-j,j},i\geq j\geq 0$.
It satisfies $L_{i,j}=L_{i-1,j}+L_{i-1,j-1}$ for $i,j\geq 1$.

\begin{prop} We have $M(\alpha)=L(\alpha)U$ where $M(\alpha),L(\alpha)$ are as above 
and where $U$ is unipotent upper-triangular with coefficients
$U={j\choose i}$ given by binomial coefficients.
\end{prop}

\noindent{\bf Proof} We have obviously $M_{i,j}=(LU)_{i,j}$
if $i=0$ or $j=0$. The remaining cases follow by induction on 
$i+j$ from the equalities
\begin{eqnarray*}
&=&\sum_k L_{i,k}U_{k,j}\\
&=&\sum_k L_{i,k}U_{k,j-1}+\sum_kL_{i,k}U_{k-1,j-1}\\
&=&\sum_k L_{i,k}U_{k,j-1}+\sum_k\tilde L_{i-1,k}U_{k,j-1}\\
&=&M_{i,j-1}+\tilde M_{i-1,j-1}\\
&=&M_{i,j-1}+M_{i-1,j}\\
&=&M_{i,j}
\end{eqnarray*}
where $\tilde M=\tilde M(1,1+\alpha_1,1+\alpha_1+\alpha_2,1+\alpha_1+\alpha_2+\alpha_3,\dots)$,
respectively $\tilde L=\tilde L(1,1+\alpha_1,1+\alpha_1+\alpha_2,\dots)$,
is obtained from $M$ by removing its first row, respectively from 
$L$ by removing its first row and column.\hfill $\Box$

\subsection{A few more identities}

The following results show other similarities between $l_{i,j}$ 
and binomial coefficients:

\begin{prop} (i)
We have for all $n,k\in\mathbb N$
the equality ${n\choose k}=l_{n,k}+2l_{n,k+1}$.

(ii) We have for all $n$ the identity
$$x^n=(-1)^n+(x+1)\sum_{k=1}^nl_{n,k}(x-1)^{k-1}.$$
\end{prop}

For proving assertion (i), 
it is enough to check the equality for all $n$ with $k=0,1$.
The general case follows from the last equality in Proposition \ref{propladd}.

The second assertion holds for $n=0$. We have now
\begin{eqnarray*}
&&(-1)^n+(x+1)\sum_{k=1}^nl_{n,k}(x-1)^{k-1}\\
&=&(-1)^n+(x+1)\sum_{k=1}^n (l_{n-1,k-1}+l_{n-1,k})(x-1)^{k-1}\\
&=&(x-1)\left((-1)^{n-1}+(x+1)\sum_{k=1}^{n-1}l_{n-1,k}(x-1)^{k-1}\right)\\
&&\quad +(-1)^{n-1}+(x+1)\sum_{k=1}^{n-1}l_{n-1,k}(x-1)^{k-1}\\
&&\quad +(-1)^n+(x+1)(-1)^{n-1}-(x-1)(-1)^{n-1}-(-1)^{n-1}\\
&=&(x-1)x^{n-1}+x^{n-1}=x^n
\end{eqnarray*}
by induction.\hfill$\Box$

\subsection{A few integer sequences}

A few integer sequences related to the numbers $l_{i,j}$
appear seemingly in \cite{OEIS} (proofs are probably easy in most cases).

Observe that the array $l_{i,j}$ with its first row removed 
appears in \cite{OEIS},
$l_{1,1},l_{2,1},l_{2,2},l_{3,1},l_{3,2},l_{3,3},l_{4,1},l_{4,2},\dots$
is A59260 of \cite{OEIS}.

The sequence $1,0,2,6,22,80,296,1106,\dots,l_{2n,n},n=0$ of 
central coefficients coincides seemingly with A72547 of \cite{OEIS}.

It is easy to show that the sequence $1,0,2,2,6,10,22,\dots,
s_n=\sum_{k=0}^nl_{n,k}$ of row-sums 
is given by $s_0=1,s_n=2(s_{n-1}+(-1)^n)=2s_{n-2}+s_{n-1}$. The 
closely related sequence
$\frac{1}{2}s_{n+1}$ coincides with A1045 of \cite{OEIS}.

$1,1,4,8,20,44,100,\dots,a_n=\sum_{k=0}^nl_{n,k}(k+1)$ coincides with 
A84219 of \cite{OEIS}.

$0,0,1,3,9,23,57,135,313,711,\dots,\sum_{k=2}^nl_{n,k}(k-1)$
coincides (up to a shift of the index) with A45883 of \cite{OEIS}.

$1,2,6,14,34,78,178,\dots,\sum_{k=0}^nl_{n,k}(2k+1)$ coincides seemingly
with A59570 of \cite{OEIS}.

There are certainly other sequences of \cite{OEIS} related to the numbers 
$l_{i,j}$. 

Interestingly, the descriptions of the above sequences are
linked to several different and apparently unrelated mathematical
areas.


\noindent{\bf Acknowledgements} I thank Bernard Parisse for a
useful discussion.


\noindent Roland BACHER, 

\noindent Univ. Grenoble Alpes, Institut Fourier, 

\noindent F-38000 Grenoble, France.
\vskip0.5cm
\noindent e-mail: Roland.Bacher@ujf-grenoble.fr

\end{document}